# [1, 2]-sets and [1, 2]-total Sets in Trees with Algorithms

A. K. Goharshady[1], M. R. Hooshmandasl[2], M. Alambardar Meybodi[3]

Department of Computer Science, Yazd University, Yazd, Iran.

The Laboratory of Quantum Information Processing, Yazd University, Yazd, Iran.

e-mail:[1] goharshady@gmail.com, [2] hooshmandasl@yazd.ac.ir, [3] alambardar_m@yahoo.com.

**Abstract**

A set $S \subseteq V$ of the graph $G = (V, E)$ is called a $[1, 2]$-set of $G$ if any vertex which is not in $S$ has at least one but no more than two neighbors in $S$. A set $S' \subseteq V$ is called a $[1, 2]$-*total* set of $G$ if any vertex of $G$, no matter in $S'$ or not, is adjacent to at least one but not more than two vertices in $S'$. In this paper we introduce a linear algorithm for finding the cardinality of the smallest $[1, 2]$-sets and $[1, 2]$-total sets of a tree and extend it to a more generalized version for $[i, j]$-sets, a generalization of $[1, 2]$-sets. This answers one of the open problems proposed in [5]. Then since not all trees have $[1, 2]$-total sets, we devise a recursive method for generating all the trees that do have such sets. This method also constructs every $[1, 2]$-total set of each tree that it generates.

**Keywords:** $[1, 2]$-sets; Dominating Sets; Trees; Total Dominating Sets.

## 1 Introduction and preliminaries

In a graph $G = (V, E)$, the (open) neighborhood of a vertex $v \in V$ is the set of all vertices adjacent to $v$ and is denoted by $N(v)$, i.e. $N(v) = \{u \in V | uv \in E\}$. The closed neighborhood of a vertex $v$ is defined $N[v] = N(v) \cup \{v\}$. If $U \subseteq V$, we define the open and close neighborhoods of $U$ as follows:

$$N(U) = \bigcup_{u \in U} N(u),$$

$$N[U] = \bigcup_{u \in U} N[u].$$

A set $S$ is called a dominating set of $G$ if $N[S] = V$, i.e. if every vertex is either in $S$ or adjacent to a vertex in $S$. The size of the smallest dominating sets of a graph $G$ is denoted by $\gamma(G)$. Any such set is called a $\gamma$-set of $G$.

There are many extensions of dominating sets such as Roman, independent, total, efficient, mixed, paired, signed and rainbow dominating sets. A discussion of some of these can be found in [11, 16]. Many researchers have contributed to this area and worked on finding bounds and algorithms to compute the cardinality of the smallest dominating sets of each kind in many families of graphs as well as the relation between domination numbers of graphs obtained by standard or well-known operations on other graphs, for example [1–3, 6–9, 13–15, 19].

A set $S \subseteq V$ is called a $[1, 2]$-set of $G$ if for each $v \in V \setminus S$ we have $1 \leq |N(v) \cap S| \leq 2$, i.e. $v$ is adjacent to at least one but not more than two vertices in $S$. Every graph has at least one $[1, 2]$-set since the set of vertices $V$ is itself a $[1, 2]$-set. The size of the smallest $[1, 2]$-sets of $G$ is denoted by $\gamma_{[1,2]}(G)$. Any such set is called a $\gamma_{[1,2]}$-set of $G$.

As a generalization of $[1, 2]$-sets, for two nonnegative integers $a$ and $b$, a set $S \subseteq V$ is an $[a, b]$-set if for each $v \in V \setminus S$ we have $a \leq |N(v) \cap S| \leq b$. A set $S' \subseteq V$ is called a total $[a, b]$-set or an $[a, b]$-total set of $G$ if for each $v \in V$ we have $a \leq |N(v) \cap S| \leq b$, i.e. if each vertex, no matter in $S'$ or not, has at least $a$ and at most $b$ neighbors in $S'$. These definitions can be found in [4, 5]. In the last section of this paper, we will mainly investigate $[1, 2]$-total sets. It is noteworthy that not every graph has such sets. This is shown in section 5.

If $u, w \in V$, then a $(u, w)$-walk of length $k$ is a sequence of (not essentially distinct) vertices $u = v_1, v_2, \ldots, v_{k+1} = w$ such that for each $1 \leq i \leq k$, $v_{i+1} \in N(v_i)$. The distance between $u$ and $w$, denoted $d(u, w)$, is the smallest $k$ such that a $(u, w)$-walk of length $k$ exists.

A graph $G = (V, E)$, where $|V| = n$, is called a tree if it is connected, has no cycles and has exactly $n - 1$ edges, i.e. $|E| = n - 1$. It is easy to show that if a simple graph satisfies two of the mentioned properties, then it must satisfy the third [17].

A tree $T$ is called a *rooted* tree, if a vertex $r \in V$ is designated as *root*. In this case for two vertices $u, v \in V$, $u$ is called the parent of $v$ if $u \in N(v)$ and $d(u, r) < d(v, r)$. Conversely, $v$ is called a child of $u$ if $v \in N(u)$ and $d(v, r) > d(u, r)$. Every vertex except the root has a unique parent [17]. A vertex $w \in V(T)$ is called a leaf if it has no children. For other terminology and concepts not explained here, refer to [16, 17].

In [5], several open problems were proposed. Some of these problems were solved in [18]. In this paper we provide a linear time algorithm for finding the cardinality of a smallest $[1, 2]$-set of a given tree, thus solving another one of the proposed problems by Chellali et. al in [5]. Then we generalize this algorithm to compute several other related sets and values, including similar values for *total* $[1, 2]$-sets. Finally, we provide a recursive method for generating all the trees that have at least one $[1, 2]$-total set.

## 2 An algorithm for computing $\gamma_{[1,2]}(T)$

We are given a tree $T$ with $n$ vertices. Our goal is to:

- Calculate $\gamma_{[1,2]}(T)$, i.e. the cardinality of a smallest $[1, 2]$-set of $T$,
- Calculate number of distinct $[1, 2]$-sets of minimal cardinality,
- Find a $[1, 2]$-set of minimal cardinality.

This section will illustrate an algorithm to reach the first two goals in linear time, i.e. $O(n)$, and the third in $O(n^2)$.

### 2.1 Definitions

We first choose an arbitrary vertex of $T$ and set it as root, so that from now on we can think of $T$ as a rooted tree. For each vertex $v$ of $T$, we define the following:

- The set $ch(v)$ consists of all children of $v$.

- $T_v$ denotes the subtree of $T$ rooted at $v$.

- $m[v]$ is the size of the smallest $[1, 2]$-set of $T_v$. This is well-defined because the set of all vertices of $T_v$ is a $[1, 2]$-set of $T_v$. This value usually appears as $\gamma_{[1,2]}(T_v)$ in literature.

- $m^+[v]$ is the size of the smallest $[1, 2]$-set of $T_v$ that contains $v$. This is well-defined because the set of all vertices of $T_v$ is a $[1, 2]$-set of $T_v$ that contains $v$.

- $U^+(v)$ is a $[1, 2]$-set of $T_v$ that contains $v$ and is of size $m^+[v]$. This set will be computed by our algorithm.

- $m^-[v]$ is the size of the smallest $[1,2]$-set of $T_v$ which does *not* contain $v$. If no such set exists we define $m^-[v]$ to be $\infty$.

- $U^-(v)$ is a $[1,2]$-set of $T_v$ that does not contain $v$ and is of size $m^-[v]$. This set will be computed by our algorithm if it exists.

- A vertex is called $(S, v)$-black, if it is included in the smallest $[1,2]$-set $S$ of $T_v$ and $(S, v)$-white otherwise. We will simply use the terms *black* and *white* when there is no risk of misunderstanding. Each such set $S$ is called a (valid) *bicoloring* of $T_v$.

It is easily inferred from the definitions above, that $m[v] = \min(m^+[v], m^-[v])$ for each vertex $v$. We are going to devise a linear time algorithm based on a bottom-up dynamic programming technique to calculate $m^-$ and $m^+$ for all vertices. To calculate these values for each vertex we may assume that similar values have already been computed for all descendants of the current vertex. This is a valid assumption since we process the vertices of $T$ in post-order.

## 2.2 Values for leaves

A vertex $v$ is called a leaf if it has no children. If $v$ is a leaf of $T$, then $T_v$ consists only of $v$ and therefore we should set $m^+[v]$ to 1 and $m^-[v]$ to $\infty$ because there's no $[1,2]$-set of $T_v$ that does not contain $v$.

## 2.3 Calculating $m^-[v]$ when $v$ is not a leaf

Since we are concerned with $m^-[v]$, we have the assumption that $v$ is white, therefore exactly one or two of its children must be black and all others white. Moreover, $v$ has no effect on how one colors the subtrees $T_u$ for $u \in ch(v)$. We consider all the possible valid bicolorings and choose the one that yields to the least number of black vertices in $T_v$ as follows:

- *Case 1:* $w$ is the only black child of $v$.

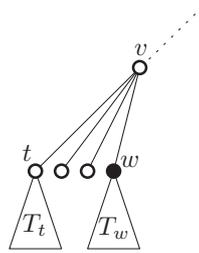

In this case we have to color $w$ black, so the minimal number of black vertices in $T_w$ will be $m^+[w]$. Similarly all other children of $v$ must be colored white, so the minimal number of black vertices in $T_t$ where $t$ is a child of $v$ other than $w$, is $m^-[t]$ and the minimal total number of black vertices in $T_v$ is:

$$c_1(v, w) := m^+[w] + \sum_{t \in ch(v) \setminus \{w\}} m^-[t]. \tag{1}$$

- *Case 2:* $u$ and $w$ are the only two black children of $v$.

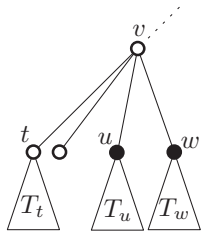

Using a similar argument as in Case 1, the minimal number of black vertices is:

$$c_2(v, u, w) := m^+[u] + m^+[w] + \sum_{t \in ch(v) \setminus \{u,w\}} m^-[t]. \qquad (2)$$

Finally, taking the minimal number of black vertices among all valid bicolorings, $m^-[v]$ can be calculated using the following formula:

$$m^-[v] = \min \left\{ \min_{w \in ch(v)} c_1(v, w), \min_{u,w \in ch(v), u \neq w} c_2(v, u, w) \right\}. \qquad (3)$$

## 2.4 Calculating $m^+[v]$ when $v$ is not a leaf

In this case, $v$ is black so all the children of $v$ already have a black neighbor and there is no restriction about their color, i.e. each child can be either white or black. The only restriction is that if $w$ is a white child of $v$, then at most one of $w$'s children can be colored black. Let $u$ and $w$ be two distinct children of $v$, then how we can bicolor $T_u$ is independent of the chosen bicoloring for $T_w$ because they are disconnected from one another. So we proceed with bicoloring each $T_w$ for $w \in ch(v)$ independently and optimally and then add up the number of black vertices in these subtrees and add 1 to the result ($v$ itself is black) to get to $m^+[v]$.

To calculate the optimal number of black vertices in $T_w$ for a specific $w \in ch(v)$, we consider the following cases:

- *Case 1:* $w$ is black:

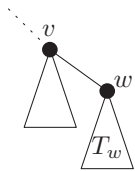

In this case the optimal number of black vertices in $T_w$ is $m^+[w]$.

- *Case 2:* $w$ is white and none of its children are black:

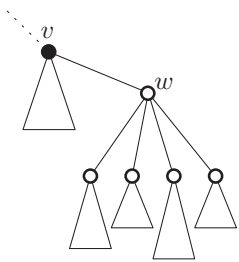

In this case we have to color all the children of $w$ white and $w$ has no effect on bicoloring the subtrees of its children (just as in 2.3), so the optimal total number of black vertices is

$$c_3(w) := \sum_{u \in ch(w)} m^-[u]. \qquad (4)$$

- *Case 3:* $w$ is white and exactly one of its children, $u$, is black:

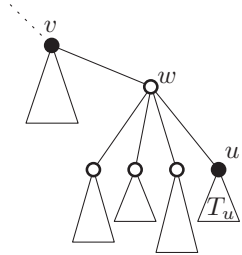

Using a similar argument as previous cases, the minimal number of black vertices in this case is:

$$c_1(w, u) = m^+[u] + \sum_{t \in ch(w) \setminus \{u\}} m^-[t]. \qquad (5)$$

We choose the optimal bicoloring in cases 1, 2 and 3 above to compute the minimal number of required black vertices in $T_w$ as follows:

$$c_4(w) := \min \left\{ m^+[w], c_3(w), \min_{u \in ch(w)} c_1(w, u) \right\}. \qquad (6)$$

So $c_4(w)$ is the minimal number of black vertices needed to bicolor $T_w$, assuming that $w$'s parent is black.

Finally $m^+[v]$ is equal to

$$1 + \sum_{w \in ch(v)} c_4(w). \qquad (7)$$

## 2.5 Example

Consider the tree below. Values of $m^+$, $m^-$ and $c_4$ can be calculated according to the following table. In the parentheses after each value the equation yielding to that value is named. Note that the vertices are numbered in post-order. Two optimal bicolorings of this tree are also provided.

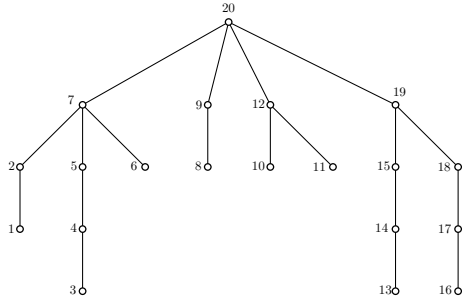

| vertex | $m^-$ | $m^+$ | $c_4$ |
|---|---|---|---|
| 1 | $\infty$ (leaf) | 1 (leaf) | 0 (leaf) |
| 2 | 1 (equation 1 ) | 1 (equation 7) | 1 (equation 5) |
| 3 | $\infty$ (leaf) | 1 (leaf) | 0 (leaf) |
| 4 | 1 (equation 1) | 1 (equation 7) | 1 (equation 5) |
| 5 | 1 (equation 1) | 2 (equation 7) | 1 (equations 4 and 5) |
| 6 | $\infty$ (leaf) | 1 (leaf) | 0 (leaf) |
| 7 | 3 (equation 1 by choosing 6 or equation 2 by choosing 2 and 6 as black children) | 3 (equation 7) | 3 (equation 5 by setting 6 as the black child) |
| 8 | $\infty$ (leaf) | 1 (leaf) | 0 (leaf) |
| 9 | 1 (equation 1) | 1 (equation 7) | 1 (equation 5) |
| 10 | $\infty$ (leaf) | 1 (leaf) | 0 (leaf) |
| 11 | $\infty$ (leaf) | 1 (leaf) | 0 (leaf) |
| 12 | 2 (equation 2) | 1 (equation 7) | 1 (same as $m^+$) |
| 13 | $\infty$ (leaf) | 1 (leaf) | 0 (leaf) |
| 14 | 1 (equation 1) | 1 (equation 7) | 1 (equation 5) |
| 15 | 1 (equation 1) | 2 (equation 7) | 1 (equations 4 and 5) |
| 16 | $\infty$ (leaf) | 1 (leaf) | 0 (leaf) |
| 17 | 1 (equation 1) | 1 (equation 7) | 1 (equation 5) |
| 18 | 1 (equation 1) | 2 (equation 7) | 1 (equations 4 and 5) |
| 19 | 3 (equation 1) | 3 (equation 7) | 2 (equation 4) |
| 20 | 8 (Several possibilities, for example by setting 12 as the only black child) | 8 (equation 7) | Not needed |

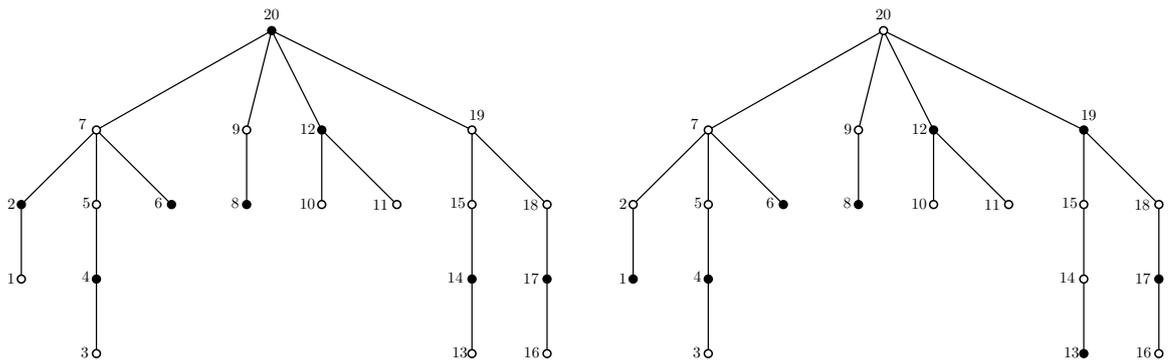

## 2.6 Time complexity

All the calculations needed to run this algorithm can be done in $O(n)$. Given a vertex $v$, $c_3(v)$ can be straightforwardly calculated in $O(|ch(v)|)$ and for each $w \in ch(v)$, we can use this equality to calculate $c_1(v, w)$ in constant time:
$$c_1(v, w) = c_3(v) + m^+[w] - m^-[w]. \tag{8}$$
So we can compute all the necessary $c_1$'s and $c_3$'s in $O(n)$.

Note that we do not need the values of all possible $c_2(v, u, w)$'s. According to equation 3, for each vertex $v$, only the minimal value of $c_2(v, u, w)$ over all possible $u$'s and $w$'s is needed. Moreover
$$c_2(v, u, w) = c_3(v) + (m^+[u] - m^-[u]) + (m^+[w] - m^-[w]), \tag{9}$$
Therefore it is sufficient to choose those vertices $u, w \in ch(v)$ that have the least $m^+[.] - m^-[.]$. This can trivially be done in $O(|ch(v)|)$ for each vertex $v$ and in $O(n)$ totally.

It is easy to show that all other actions of this algorithm take no more than linear time. So the algorithm runs in $O(n)$ overall.

Note that some of the equations above are not valid in presence of infinite values and need special treatment in implementation. We did not mention this in this section for the sake of simplicity. Consult the pseudo-code in the appendix about this matter. Similar techniques should be applied to implement other algorithms in this text as well. These will neither change the rationale behind the algorithm nor its runtime.

## 2.7 Finding a $\gamma_{[1,2]}$-set and number of such sets

To find a smallest $[1, 2]$-set of $T$, we first find a smallest set between all those $[1, 2]$-sets of $T$ that contain the root, then we find a smallest set between those that do not contain the root. The smaller of these two is the set we are looking for. In other words, we compute $U^+(root)$ and $U^-(root)$ and take the one with smaller cardinality (any one is ok if they are of equal size).

We compute $U^-(v)$ and $U^+(v)$ for some vertices $v$, including the root, using a recursive method that is described below. In using this method, we assume that we have saved all the necessary information after running the algorithm that was explained above.

If $m^-[v] = \infty$, for example when $v$ is a leaf, then $U^-(v)$ doesn't exist and the algorithm terminates.

Otherwise, we check our calculation of equation 3 for $v$ to see which children of $v$ were supposed to be colored black. We proceed by finding $U^+(w)$ for each black child $w$ of $v$ and $U^-(t)$ for each white child $t$. Since $v \notin U^-(v)$, $U^-(v)$ is the union of these sets.

Going to $U^+(v)$, we start by setting $\{v\}$ as our estimate of $U^+(v)$. According to the idea in section 2.4, we have to find a suitable bicoloring for each $T_w$ where $w$ is a child of $v$. For each such $w$, we check equation 6 for the minimal term. If $m^+[w]$ was minimal, then we add $U^+(w)$ to our estimate of $U^+(v)$. Otherwise if $c_3(w)$ was minimal, we add $U^-(u)$ for all $u \in ch(w)$ to our set. Otherwise if $u$ is supposed to be the only black child of $w$, we add $U^+(u)$ and $U^-(t)$ for all children $t$ of $w$ except $u$ to our set.

This process is $O(n^2)$ because each set $U^+(v)$ or $U^-(v)$ can a have a size of at most $n$. Its correctness can be proven according to the definitions, in the same manner as of the previous algorithm. Note that this algorithm can be edited to find *all* smallest $[1, 2]$-sets of $T$, taking $O(n^2)$ time per answer. We will later show that number of different answers, i.e. different smallest $[1, 2]$-sets, in a given tree may not be bounded by a polynomial in terms of $n$ – number of vertices.

Same idea, which is sometimes called back-substitution, can be applied to the problem of finding the number of $\gamma_{[1,2]}$-sets. Start by defining $\nu^+(v)$ as the number of different choices that exist for $U^+(v)$. Define $\nu^-(v)$ similarly. Then these numbers can be calculated using the same recursive method, by applying the addition principle to all the paths that lead to a minimal value in equations 3 and 6 and using the technique in section 2.6 and applying the multiplication principle according to equations 1, 2, 4, 5, and 7. This algorithm is linear.

## 2.8 Cardinality of a smallest $[a, b]$-set

The same algorithm can be slightly modified to work for general $[a, b]$-sets. For each vertex $v$ of $T$ we define the following:

- $m_{[a,b]}[v]$ is the size of the smallest $[a, b]$-set of $T_v$.

- $m^+_{[a,b]}[v]$ is the size of the smallest $[a, b]$-set of $T_v$ that contains $v$.

- $m^-_{[a,b]}[v]$ is the size of the smallest $[a, b]$-set of $T_v$ that does *not* contain $v$.

Again, it is obvious from the definition that $m_{[a,b]}[v] = \min\left\{m^+_{[a,b]}[v], m^-_{[a,b]}[v]\right\}$. The initial values, i.e. values for leaves, can be computed in exactly the same manner as the algorithm for $[1, 2]$-sets. We proceed with computing $m^-_{[a,b]}[v]$ for those vertices $v$ of $T$ that are not leaves.

Since $v$ is white, number of its black children must be in range $[a, b]$. Let $A$ denote the set of black children of $v$. Then the minimal number of black vertices in $T_v$ is:

$$e_1(A) := \sum_{u \in A} m^+_{[a,b]}[u] + \sum_{u \in ch(v) \setminus A} m^-_{[a,b]}[u] \tag{10}$$

$$= \sum_{u \in ch(v)} m^-_{[a,b]}[u] + \sum_{u \in A} \left(m^+_{[a,b]}[u] - m^-_{[a,b]}[u]\right). \tag{11}$$

We can calculate the first sum for all $v$'s in $O(n)$ and the minimal value for the second sum can be found according to the technique in section 2.6. Since we have to find between $a$ and $b$ optimal black vertices, and not only one or two, that minimize $m^+_{[a,b]}[.] - m^-_{[a,b]}[.]$, we can use a heap for storing these values. This will yield to a total runtime of $O(n \log b)$.

To achieve this runtime, note that we can insert values of $m^+_{[a,b]}[u] - m^-_{[a,b]}[u]$ for all $u$'s into a max-heap and then extract the required sum from this heap, but whenever the size of our heap exceeded $b$, we can remove a maximal element from the heap, thus keeping at most $b$ items in the heap at each moment. Note that we can always store the sum of all elements in our heap separately and will therefore not need to visit all the items in the heap to find this sum. Inserting and deleting one item in a heap of size $b$ takes $O(\lg b)$ time [12]. For each vertex $u$, at most one item is inserted to one such heap and possibly removed from it, so the total runtime of our heap operations is no more than $O(n \lg b)$. Other parts of the algorithm are similar to previous ones and are linear.

To calculate $m^+_{[a,b]}[v]$, for a non-leaf $v$, we can use the same approach as when we were computing $m^+[v]$ in section 2.4. In doing so, one should mind that since $v$ is a black vertex, if $u$ is a white child of $v$, then number of black children of $u$ must be between $a - 1$ and $b - 1$ inclusively. This can again be done in $O(n \log b)$ in a similar manner.

## 3 Algorithm for total sets

In this section we tweak the algorithm in previous section to find the cardinality of a smallest total $[1, 2]$-set of a given tree $T$. Just as in last section, one can find the number of smallest total $[1, 2]$-sets and the sets themselves using back-substitution. Moreover, this extension of the algorithm is not restricted to $[1, 2]$-sets and can be applied to find both cardinality and number of smallest $[a, b]$-sets, for any nonnegative integers $a \leq b$, in a tree.

## 3.1 New definitions

We keep the definitions of the previous section. Let $k \in \{1, 2\}$. We define the following for each vertex $v$ of $T$:

- $m_k^-[v]$ is defined to be the cardinality of a smallest total $[1,2]$-set of $T_v$ that does not contain $v$ and contains at most $k$ of $v$'s children. If there is no such set, we let this number be $\infty$.

- $m_k^+[v]$ is defined to be the cardinality of a smallest total $[1,2]$-set of $T_v$ that contains $v$ and at most $k$ of $v$'s children. If there is no such set, we let this number be $\infty$.

- If $S$ is a total $[1,2]$-set of a graph $G$ and $u \in S$, we say that $u$ is $(S, G)$-total black. If $u \notin S$, we say that $u$ is $(S, G)$-total white. $S$ is called a total bicoloring of $G$. We simply use the words "black", "white" and "bicoloring" when there is no risk of misunderstanding.

- If $S$ is a total bicoloring of $G$ and $u \in S$ and $v \in N(u)$, we say that $u$ covers $v$ in $S$. We may omit $S$ and just say that $u$ covers $v$ when there is no risk of ambiguity.

We use the same bottom-up approach. According to these definitions, the cardinality of a smallest total $[1,2]$-set of $T$ is:

$$\min \left\{ m_2^+[root], m_2^-[root] \right\}. \tag{12}$$

Note that for each vertex $v$ of $T$, $m_1^-[v] \geq m_2^-[v]$ and $m_1^+[v] \geq m_2^+[v]$. If the value of (12) tends to be $\infty$, then the tree has no total $[1,2]$-sets. We will further investigate the set of trees that have at least one total $[1,2]$-set in next sections. The algorithm has the same complexity as the previous one. This can be proven using the same statement as of the proof in 2.6.

## 3.2 Values for leaves

Suppose $v$ is a leaf of $T$, So $T_v$ consists of $v$ only and has no total $[1,2]$-sets. Therefore $m_1^-[v] = m_2^-[v] = m_1^+[v] = m_2^+[v] = \infty$.

## 3.3 Calculating $m_k^-[v]$ when $v$ is not a leaf

In this case $v$ is white. At least one and at most $k$ of its children must be colored black and all of them must be covered by their own children, since $v$ does not cover them. Let $A$ be the set of black children of $v$, then the minimal number of black vertices is:

$$\sum_{w \in A} m_2^+[w] + \sum_{w \in ch(v) \setminus A} m_2^-[w],$$

so

$$m_k^-[v] = \min_{A \subseteq ch(v), 1 \leq |A| \leq k} \sum_{w \in A} m_2^+[w] + \sum_{w \in ch(v) \setminus A} m_2^-[w].$$

**Remark:** These summations are assumed to be zero when the summand is changing over the empty set.

## 3.4 Calculating $m_k^+[v]$ when $v$ is not a leaf

$v$ is black, so it covers all its children and if $w \in ch(v)$, then $w$ can have at most one black child. $v$ must have at least one and at most $k$ black children. We proceed by optimally bicoloring $T_w$ for each vertex $w \in ch(v)$.

- *Case 1:* $w$ is white. Since $w$ is covered by $v$, it has either one or no black children. If it has one black child, then the optimal bicoloring takes $m_1^-[w]$ black vertices, otherwise since all the children of $w$ are white and are not covered by $w$, the minimal number of needed black vertices is $\sum_{u \in ch(w)} m_2^-[u]$. So the total minimal number of black vertices needed to color $T_w$ in this case is:

$$d_1(w) := \min \left\{ m_1^-[w], \sum_{u \in ch(w)} m_2^-[u] \right\}.$$

- *Case 2:* $w$ is black. If $w$ is covered by at least one of its children, we will need $m_1^+[w]$ black vertices, otherwise all its children are white and are covered by $w$, so we need $1 + \sum_{u \in ch(w)} m_1^-[u]$ black vertices including $w$. So the minimal number of black vertices in this case is:

$$d_2(w) := \min \left\{ m_1^+[w], 1 + \sum_{u \in ch(w)} m_1^-[u] \right\}.$$

Finally if $A$ is supposed to be the set of black children of $v$,

$$m_k^+[v] = \min_{A \subseteq ch(v), 1 \leq |A| \leq k} \sum_{w \in A} d_2(w) + \sum_{w \in ch(v) \setminus A} d_1(w).$$

## 3.5 Cardinality of a smallest total $[a, b]$-set

To find the cardinality of a smallest total $[a, b]$-set of a given tree $T$, Let $k \in \{b, b-1\}$ it is sufficient to redefine the variables as follows:

- $m_k^-[v]$ is defined to be the cardinality of a smallest total $[a, b]$-set of $T_v$ that does not contain $v$ and contains at most $k$ of $v$'s children. If there is no such set, we let this number be $\infty$.

- $m_k^+[v]$ is defined to be the cardinality of a smallest total $[a, b]$-set of $T_v$ that contains $v$ and at most $k$ of $v$'s children. If there is no such set, we let this number be $\infty$.

We have the same values for leaves and the other equations change like this:

$$m_k^-[v] = \min_{A \subseteq ch(v), a \leq |A| \leq k} \sum_{w \in A} m_b^+[w] + \sum_{w \in ch(v) \setminus A} m_b^-[w],$$

$$d_1'(w) := \min \left\{ m_{b-1}^-[w], \sum_{u \in ch(w)} m_b^-[u] \right\},$$

$$d_2'(w) := \min \left\{ m_{b-1}^+[w], 1 + \sum_{u \in ch(w)} m_b^-[u] \right\},$$

$$m_k^+[v] = \min_{A \subseteq ch(v), a \leq |A| \leq k} \sum_{w \in A} d_2'(w) + \sum_{w \in ch(v) \setminus A} d_1'(w).$$

These can be proved in exactly the same manner as in total $[1, 2]$-sets.

# 4 On the number of $\gamma_{[1,2]}$-sets in a tree

In 2.7, we discussed an algorithm that was able to find *all* $\gamma_{[1,2]}$-sets of a tree $T$ in $O(n^2 \times n_\gamma(T))$ where $n$ is the number of vertices of $T$ and $n_\gamma(T)$ denotes the number of distinct smallest $[1,2]$-sets, i.e. $\gamma_{[1,2]}$-sets, of $T$. The natural problem that arises from that discussion is whether there exists a polynomial algorithm to find *all* such sets. In this section we show that $n_\gamma(T)$ may not be bounded by a polynomial of $n$ and so, the answer to the previous question is "No".

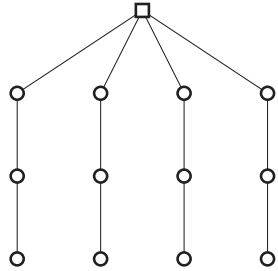

We are going to consider a special family of trees like the figure above. These trees consist of a root vertex, depicted as a square, which has $k$ children and if a vertex $u$ is a child of the root, then $T_u$ is isomorphic to $P_3$. For simplicity, we call each of these $P_3$'s a "branch" of the tree. One can easily find a $[1,2]$-set of size $k+1$ in each such tree. It is sufficient to color the root and the middle vertex of each branch black.

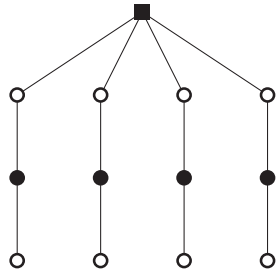

It can also be observed that, since the bicolorings of different branches are independent, one can color the root black and then in each branch set the color of either the middle vertex or the leaf to black. Any such bicoloring gives a $[1,2]$-set of $k+1$ black vertices. Therefore, we have at least $2^k$ such sets.

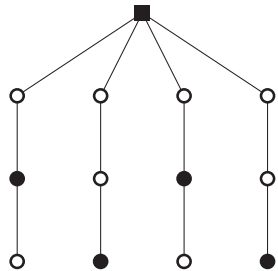

Now we prove that there exists no $[1,2]$-set of size $k$ or less. We consider two cases. If the root is black, then we must have at least one black vertex in each branch for a total of at least $k+1$ blacks. If the root is white, then one or two of its children must be black and we will still need one more black vertex in each of the branches, no matter the color of the upper-most vertex of that branch is white or black, again forcing us to use at least $k+1$ blacks.

To sum up, such graphs have $3k+1$ vertices and at least $2^k$ distinct $\gamma_{[1,2]}$-sets, so the number of such sets cannot be bounded by a polynomial of the number of vertices in a tree.

## 5 $[1,2]$-total trees

A tree $T$ is called a $[1,2]$-total tree if it has at least one total $[1,2]$-set. Not all trees are $[1,2]$-total, the graph below is an example of a tree that lacks total $[1,2]$-sets:

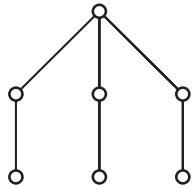

In [10], Galvas and Shultz construct a recursive method for identifying the family of all efficient open domination trees. Similarly, we devise another recursive method that generates all $[1,2]$-total trees. Our method also provides a valid bicoloring for each tree it generates, i.e. along with each tree, the algorithm constructs a total $[1,2]$-set for that tree as well. Let's denote the set of all $[1,2]$-total trees as $\Upsilon$.

The method constructs $\Upsilon$ as follows:

- $P_2$ is in $\Upsilon$. A total $[1,2]$-set for $P_2$ is the set consisting of its both vertices.

- If $T \in \Upsilon$ and $S$ is a total $[1,2]$-set of $T$ constructed by this algorithm, then the following are also in $\Upsilon$:

  1. A tree $T_1$ that is obtained from adding a new white vertex to $T$ and joining the new vertex with a vertex in $S$. $S$ is a total $[1,2]$-set for each such $T_1$.

  2. If $v$ is a black vertex of $T$, i.e. $v \in S$, and has exactly one black neighbor, then a new tree $T_2$ can be obtained by adding a new black vertex $w$ to $T$ and joining $v$ and $w$. A total $[1,2]$-set for the new tree is $S \cup \{w\}$.

  3. Let $u \in V(T) \setminus S$. If $u$ has only one black neighbor then we can add two or three new vertices, $v$, $w$ and $x$ to our tree according to the figure below:

     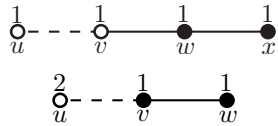

     In this figure, the number of black neighbors of each vertex is written next to it. Black vertices are added to $S$ to form a total $[1,2]$-set for the new graph.

  4. Let $u \in V(T) \setminus S$ be a vertex having two black neighbors. We can add three new vertices according to this figure:

     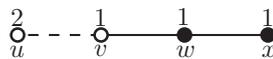

- Any tree in $\Upsilon$ can be constructed using a finite number of applications of the rules above.

Along with each tree we construct, we provide a total $[1,2]$-set of it as well, so it is obvious that every tree made by this algorithm is indeed in $\Upsilon$. Now we prove that the converse holds as well, i.e. the method above generates every tree in $\Upsilon$. Moreover, we prove that the method described above finds all valid bicolorings for each of the trees in $\Upsilon$.

There is no $[1,2]$-total tree with less than two vertices. The only $[1,2]$-total tree with two vertices is $P_2$ and the only bicoloring of $P_2$ is the one consisting of both vertices. Let $T$ be a $[1,2]$-total tree with $n \geq 3$ vertices, and let $S$ be an arbitrary total $[1,2]$-set of $T$. Using induction, we show that our algorithm generates $T$. Assume that the algorithm generates any $[1,2]$-total tree having less than $n$ vertices.

If $T$ has a leaf that is not in $S$, then construct the tree $T'$ by removing this leaf from $T$. $S$ is a total $[1,2]$-set of $T'$ as well, and $T'$ has less than $n$ vertices, so our algorithm constructs $T'$. $T$ can be obtained from $T'$ by applying rule number 1 above.

Assume that $S$ contains all the leaves of $T$. Let $P$ be a longest path of $T$ and let $u$, $v$ and $w$ be its first three vertices.

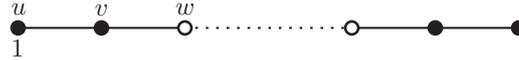

$u$ is assumed to be in $S$ and $v$ must be in $S$ because it is $u$'s only neighbor. If $v$ has a black neighbor other than $u$, then we can remove $u$ from the tree, construct the rest of the tree by induction hypothesis and then add $u$ according to rule number 2. So suppose that all other neighbors of $v$, including $w$, are white. If $v$ is the only black neighbor of $w$, we can delete the vertices $u$, $v$ and $w$, construct the rest of the tree and then add them according to rules 3 or 4. Note that in this case $w$ can have at most two neighbors, including $v$. Otherwise, i.e. if $w$ has two black neighbors, we can omit $u$ and $v$, make the rest of the tree and then add them back according to part 2 of rule number 3. This concludes the proof.

# References

[1] R. B. Allan and R. Laskar. On domination and independent domination numbers of a graph. *Discrete Mathematics*, 23(2):73–76, 1978.

[2] B. Bresar, M. A. Henning, and D. F. Rall. Rainbow domination in graphs. *Taiwanese Journal of Mathematics*, 12(1):213–225, 2008.

[3] G. J. Chang, B. S. Panda, and D. Pradhan. Complexity of distance paired-domination problem in graphs. *Theoretical Computer Science*, 459(0):89 – 99, 2012.

[4] M. Chellali, O. Favaron, T. W. Haynes, and S. T. Hedetniemi. Independent $[1,k]$-sets in graph. *Australian Journal of Combinatorics*, 59(1):144–156, 2014.

[5] M. Chellali, T. W. Haynes, S. T. Hedetniemi, and A. McRae. $[1,2]$-sets in graphs. *Discrete Applied Mathematics*, 161(18):2885–2893, December 2013.

[6] L. Chen, C. Lu, and Z. Zeng. Hardness results and approximation algorithms for (weighted) paired-domination in graphs. *Theoretical Computer Science*, 410(4749):5063 – 5071, 2009.

[7] E. J. Cockayne, Dawes R. M., and S. T. Hedetniemi. Total domination in graphs. *Networks*, 10(3):211–219, Autumn (Fall) 1980.

[8] W. Duckworth and N. C. Wormald. Minimum independent dominating sets of random cubic graphs. *Random Structures & Algorithms*, 21(2):147–161, September 2002.

[9] M. R. Fellows and Hoover M. N. Perfect domination. *Australasian Journal of Combinatorics*, 3:141–150, 1991.

[10] H. Gavlas and K. Schultz. Efficient open domination. *Electronic Notes in Discrete Mathematics*, 11:681–691, July 2002.

[11] M. A. Henning and A. Yeo. *Total Domination in Graphs*. Springer Monographs in Mathematics, 2013.

[12] E. Horowitz and S. Sahni. *Fundamentals of Data Structures*, 357–359. Computer Science Press, 1982.

[13] R. W. Irving. On approximating the minimum independent dominating set. *Information Processing Letters*, 37(4):197–200, February 1991.

[14] J. K. Lan and G. J. Chang. On the mixed domination problem in graphs. *Theoretical Computer Science*, 476(0):84 – 93, 2013.

[15] M. Mollard. The domination number of cartesian product of two directed paths. *Journal of Combinatorial Optimization*, 27(1):144–151, January 2014.

[16] T. W. Haynes, S. T. Hedetniemi, and P. Slater. *Fundamentals of Domination in Graphs*. CRC Press, 1998.

[17] D. B. West. *Introduction to Graph Theory*. Pearson, 2nd edition, September 2000.

[18] X. Yang and B. Wu. $[1,2]$-domination in graphs. *Discrete Applied Mathematics*, 175:79–86, October 2014.

[19] Y. Zhao, L. Kang, and M. Y. Sohn. The algorithmic complexity of mixed domination in graphs. *Theoretical Computer Science*, 412(22):2387 – 2392, 2011.